\documentclass[11pt]{amsart}

\usepackage{amsmath}
\usepackage{amssymb}
\usepackage{enumerate}
\usepackage{euscript}
\usepackage{amscd}
\usepackage[all]{xy}

\newtheorem{thm}{Theorem}[section]
\newtheorem{lem}[thm]{Lemma}
\newtheorem{prop}[thm]{Proposition}

\theoremstyle{definition}
\newtheorem{defn}[thm]{Definition}
\newtheorem{rem}[thm]{Remark}

\newtheorem{asmp}[thm]{Assumption}

\numberwithin{equation}{section}
\theoremstyle{remark}

\newcommand{\bba}{{\mathbb A}}

\newcommand{\bbc}{{\mathbb C}}

\newcommand{\bbq}{{\mathbb Q}}
\newcommand{\bbr}{{\mathbb R}}

\newcommand{\bbz}{{\mathbb Z}}

\newcommand{\gA}{{\mathfrak A}}

\newcommand{\gB}{{\mathfrak B}}
\newcommand{\gc}{{\mathfrak c}}
\newcommand{\gC}{{\mathfrak C}}

\newcommand{\gG}{{\mathfrak G}}

\newcommand{\gM}{{\mathfrak M}}

\font\tenscr=rsfs10 

\newcommand{\sC}{\hbox{\tenscr C}}

\newcommand{\sE}{\hbox{\tenscr E}}

\newcommand{\sS}{\hbox{\tenscr S}}

\newcommand{\esA}{{\EuScript{A}}}
\newcommand{\esB}{{\EuScript{B}}}
\newcommand{\esC}{{\EuScript{C}}}

\newcommand{\cK}{{\mathcal K}}

\newcommand{\co}{{\mathcal O}}

\newcommand{\cZ}{{\mathcal Z}}

\newcommand{\rd}{{\operatorname{rd}}}
\newcommand{\ird}{{\operatorname{ird}}}

\newcommand{\aut}{{\operatorname{Aut}}\,}

\newcommand{\aff}{{\operatorname {Aff}}}

\newcommand{\re}{{\operatorname {Re}}}

\newcommand{\gl}{{\operatorname{GL}}}

\newcommand{\sst}{{\operatorname{ss}}}

\newcommand{\sep}{{\operatorname{sep}}}

\newcommand{\br}{{\operatorname{B}}}

\newcommand{\st}{{\operatorname{St}}}
\newcommand{\diag}{{\operatorname{diag}}}
\newcommand{\cl}{{\operatorname{Cl}}}

\newcommand{\bs}{Schwartz--Bruhat function}

\newcommand{\pv}{prehomogeneous vector space}

\newcommand{\res}{\operatorname{Res}}

\newcommand{\A}{\bba}
\newcommand{\Z}{\bbz}
\newcommand{\Q}{\bbq}
\newcommand{\R}{\bbr}
\newcommand{\C}{\bbc}

\newcommand{\ma}{\bba^{\times}}

\newcommand{\fa}{\mathbb A_{\mathrm f}}
\newcommand{\mfa}{\mathbb A_{\mathrm f}^\times}
\newcommand{\mk}{k^{\times}}

\newcommand{\md}{d^{\times}}
\newcommand{\gMf}{\gM_{\rm f}}

\newcommand{\fin}{{\rm f}}
\newcommand{\pr}{\mathrm{pr}}

\newcommand{\twtw}[4]
{\begin{pmatrix}{#1}&{#2}\\{#3}&{#4}\\\end{pmatrix}}
\newcommand{\stwtw}[4]
{\scriptsize\begin{pmatrix}{#1}&{#2}\\{#3}&{#4}\\\end{pmatrix}}

\newcommand{\aaa}{\mathfrak a}
\newcommand{\ccc}{\mathfrak c}
\newcommand{\gaa}{G_\aaa}
\newcommand{\vaa}{V_\aaa}
\newcommand{\bac}{B_{\aaa,\ccc}}
\newcommand{\wac}{W_{\aaa,\ccc}}

\renewcommand{\sC}{\hbox{\tenscr C}\, }

\begin{document}

\title[distributions of cubic algebras]
{Distributions of discriminants of cubic algebras II}

\author[Takashi Taniguchi]{Takashi Taniguchi}
\address{Department of Mathematical Sciences, University of Tokyo}
\email{tani@ms.u-tokyo.ac.jp}
\thanks{Financial support is provided by
Research Fellowships for Young Scientists
of Japan Society for the Promotion of Science.}
\date{\today}

\begin{abstract}
Let $k$ be a number field and $\co$ the ring of integers.
In the previous paper \cite{ddca} we study the Dirichlet series
counting discriminants of cubic algebras of $\co$
and derive some density theorems on distributions of the discriminants
by using the theory of zeta functions of prehomogeneous vector spaces.
In this paper we consider these objects under imposing
finite number of splitting conditions at non-archimedean places.
Especially the explicit formulae of residues at $s=1$ and $5/6$
under the conditions are given.
\end{abstract}


\maketitle


\section{Introduction}\label{sec:introduction}


Let $k$ be a number field and $\co$ the ring of integers.
For a place $v$ of $k$ let $k_v$ be
the completion of $k$ at $v$.
Let $T$ be a finite set of places.
Take a separable cubic algebra $L_v$ of $k_v$ for each $v\in T$
and let $L_T=(L_v)_{v\in T}$ the $T$-tuple.
We call an $\co$-algebra a {\em cubic algebra}
if it is locally free of rank $3$ as an $\co$-module.
We denote by $\sC(\co)$ the set of isomorphism
classes of cubic algebras of $\co$.
We let
\[
\sC(\co,L_T)^\ird:=
\left\{R\in\sC(\co)\ \ \vrule
	\begin{array}{l}
	\text{$F=R\otimes_\co k$ is a cubic field extension of $k$, and}\\
	\text{$F\otimes_k k_v\cong L_v$ for all $v\in T$.}
	\end{array}
\right\}.
\]
We define
\begin{align*}
\vartheta_{L_T}^\ird(s)
&:=	\sum_{R\in \sC(\co,L_T)^\ird}
		\frac{{}^\#(\aut(R))^{-1}}{N(\Delta_{R/\co})^s},\\
h_{L_T}(X)
&:={}^\#\{R\in\sC(\co,L_T)^\ird\mid N(\Delta_{R/\co})<X\}.
\end{align*}
Here we denote by $N(\Delta_{R/\co})$ the ideal norm
of the relative discriminant of $R/\co$
and by $\aut(R)$ the group of automorphisms
of $R$ as an $\co$-algebra.
The primary purpose of this paper is to prove the following.
Let $n=[k:\Q]$.
\begin{thm}\label{thm:mainthm_total}
There exist constants $\gA_{L_T}$ and $\gB_{L_T}$
described explicitly such that;
\begin{enumerate}[{\rm (1)}]
\item
$\vartheta_{L_T}^\ird(s)$
has meromorphic continuation to the whole complex plane
which is holomorphic for ${\rm Re}(s)>1/2$ except for simple poles
at $s=1$ and $5/6$ with residues $\gA_{L_T}$ and $\gB_{L_T}$,
respectively, and
\item for any $\varepsilon>0$,
\[
h_{L_T}(X)
	=\gA_{L_T}X+(5/6)^{-1}\gB_{L_T}X^{5/6}
		+O(X^{\frac{5n-1}{5n+1}+\varepsilon})
\qquad (X\to\infty).
\]
\end{enumerate}
\end{thm}
Note that the $X^{5/6}$-term in the formula is relevant only when $n=1,2$.
We give the formulae of $\gA_{L_T}$ and $\gB_{L_T}$.
Let $\zeta_k(s)$ and $\Delta_k$ be the Dedekind zeta function
and the absolute discriminant of $k$, respectively.
We denote by $\gM_\infty$, $\gM_\R$, $\gM_\C$ and $\gMf$
the set of all infinite places, real places, complex places
and finite places, respectively.
We put $r_1={}^\#\gM_\R$ and $r_2={}^\#\gM_\C$.
For $v\in\gMf$, let $q_v$ be the order of the residue field of $k_v$.
We put $\theta_{L_v}={}^\#(\aut_{k_v\text{-algebra}}(L_v))$.
For a non-archimedean local field $K$ with the order of residue field $q$,
we define its local zeta function by $\zeta_K(s)=(1-q^{-s})^{-1}$.
The cubic algebra $L_v$ is in general a product of local fields.
We define $\zeta_{L_v}(s)$ as the product
of the zeta functions of those fields.
The relative discriminant $\Delta_{L_v/k_v}$
is also defined as the product of relative discriminants
of those local fields.
We denote by $\Delta_{L_v}$ its norm.
We put $i_\infty(L_T)=	{}^\#\{v\in\gM_\R\mid L_v=\R^3\}$.
We give the value in case of $T\supset\gM_\infty$.
The general case 
is easily obtained from this by taking a suitable summation.

\begin{thm}\label{thm:residueformula}
When $T\supset\gM_\infty$,
the constants $\gA_{L_T}$ and $\gB_{L_T}$ are given by
\begin{align*}
\gA_{L_T}
&=	\frac{\res_{s=1}\zeta_k(s)\cdot \zeta_k(2)}{2^{r_1+r_2+1}3^{i_\infty(L_T)+r_2}}
\prod_{v\in T\cap\gMf}
\alpha_v(L_v),\\
\gB_{L_T}
&=	\frac{\res_{s=1}\zeta_k(s)\cdot{\zeta_k(1/3)}}
{{6\Delta_k^{1/2}(\sqrt3)^{r_2-i_\infty(L_T)}}}
	\left(\frac{3\Gamma(1/3)^3}{2\pi}\right)^n
\prod_{v\in T\cap\gMf}
\beta_v(L_v),
\end{align*}
where
\begin{align*}
\alpha_v(L_v)
&=	
	\frac{(1-q_v^{-1})(1-q_v^{-2})}{(1-q_v^{-4})(1-q_v^{-5})}
	\cdot	\theta_{L_v}^{-1}\Delta_{L_v}^{-1}
	\cdot	\frac{\zeta_{L_v}(2)}{\zeta_{L_v}(4)},\\
\beta_v(L_v)
&=	
	\frac{(1-q_v^{-1/3})(1-q_v^{-1})}{(1-q_v^{-10/3})(1-q_v^{-4})}
	\cdot	\theta_{L_v}^{-1}\Delta_{L_v}^{-1}
	\cdot	\frac{\zeta_{L_v}(1/3)\zeta_{L_v}(5/3)}
			{\zeta_{L_v}(2/3)\zeta_{L_v}(10/3)}.
\end{align*}
\end{thm}

Let $v\in\gMf$. We see by computation that
$\sum_{L_v}\alpha_v(L_v)=\sum_{L_v}\beta_v(L_v)=1$ where
$L_v$ runs through all the separable cubic algebras of $k_v$.
Hence $\alpha_v(L_v)$ and $\beta_v(L_v)$ give the proportion
of the contributions of cubic algebras with local splitting type $L_v$.
The computation of $\alpha_v(L_v)$ is reduced to the determination
of certain orbital volume in a $v$-adic vector space.
The meaning of $\beta_v(L_v)$ is more subtle
and the computation requires a careful local theory.

In this paper we also try to improve Theorem \ref{thm:mainthm_total}
under segregating the cubic algebras via the {\em Steinitz class}.
The Steinitz class takes value in the ideal class group,
and it is known that finitely generated locally free modules
over a Dedekind domain
is completely classified by the rank and the Steinitz class
(cf.\ \cite{milnorb}.)
Let $\cl(k)$ be the ideal class group of $k$.
For $R\in\sC(\co)$, let $\st(R)\in\cl(k)$ be its Steinitz class.
This is the class of a fractional ideal isomorphic to $\bigwedge^3R$
as an $\co$-module.
For $\mathfrak a\in\cl(k)$, we define
\begin{align*}
\vartheta_{L_T}^\ird(\mathfrak a,s)
&:=	\sum_{R\in\sC(\co,L_T)^\ird,\ \st(R)=\mathfrak a}
		\frac{{}^\#(\aut(R))^{-1}}{N(\Delta_{R/\co})^s},\\
h_{L_T}(\mathfrak a, X)
&:={}^\#\{R\in\sC(\co,L_T)^\ird\mid \st(R)=\mathfrak a,
	\ N(\Delta_{R/\co})<X\}.
\end{align*}
Then obviously
$\vartheta_{L_T}^\ird(s)=
\sum_{\mathfrak a\in\cl(k)}\vartheta_{L_T}^\ird(\mathfrak a,s)
$
and 
$h_{L_T}(X)=\sum_{\mathfrak a\in\cl(k)}h_{L_T}(\mathfrak a, X)$.
Let $h_k$ be the class number of $k$ and
$h_k^{(3)}$ the number of $3$-torsions of $\cl(k)$
(which is a power of $3$.)
Also for $\mathfrak a\in \cl(k)$, we put $\tau(\mathfrak a)=1$
if there exists
$\mathfrak b\in \cl(k)$ such that $\mathfrak a=\mathfrak b^3$
and $\tau(\mathfrak a)=0$ otherwise.
In this paper we also prove the following.
\begin{thm}\label{thm:mainthm_stnz}
\begin{enumerate}[{\rm (1)}]
\item
The Dirichlet series
$\vartheta_{L_T}^\ird(\mathfrak a,s)$
has a meromorphic continuation to the whole complex plane
which is holomorphic for ${\rm Re}(s)>1/2$ except for a simple poles
at $s=1$ with the residue $\gA_{L_T}/h_k$ and
a possible simple pole at $s=5/6$ with the residue
$\tau(\mathfrak a)\gB_{L_T}h_k^{(3)}/h_k$,
respectively.
\item Assume $L_T$ is chosen so that at least one of $L_v$ is a field.
Then $\vartheta_{L_T}^\ird(\mathfrak a,s)$ is holomorphic in the whole complex
plane except for $s=1$ and $5/6$. Also
for any $\varepsilon>0$,
\[
h_{L_T}(\mathfrak a,X)
	=\frac{\gA_{L_T}}{h_k}X
	+\tau(\mathfrak a)\frac{h_k^{(3)}}{h_k}\cdot
		\frac{\gB_{L_T}}{5/6}X^{5/6}
		+O(X^{\frac{4n-1}{4n+1}+\varepsilon})
\qquad (X\to\infty).
\]
\end{enumerate}
\end{thm}
The statement (1) is an improvement of
Theorem \ref{thm:mainthm_total} (1).
In case $[k:\Q]=2$ and $3\mid h_k$,
the formula in (2) implies that the distribution of the Steinitz classes
of the elements of $\sC(\co,L_T)^\ird$
is not uniform, and the irregularity is
reflected in the $X^{5/6}$-term.
Although it is likely that such a formula exists for unconditional $L_T$ also,
to remove the condition on $L_T$ is highly non-trivial.

In the previous paper \cite{ddca} we prove these theorems for $T=\gM_\infty$.
(The original Shintani's theorem \cite[Theorem 4]{shintanib}
is for $k=\Q$, $T=\{\infty\}$.)
As in \cite{ddca}
we approach these theorems by using
{\em Sato-Shintani's zeta function} \cite{sash}
for {\em the space of binary cubic forms}. 

The study of class numbers of integral binary cubic forms over $\Z$
was initiated by G. Eisenstein and developed by many mathematicians
including C. Hermite, H. Davenport and T. Shintani.
Via Delone-Faddeev's correspondence \cite{defa} the $\gl_2(\Z)$-orbits
of integral binary cubic forms corresponds bijectively
to the set of cubic rings, thus useful to investigate
cubic fields or their orders.
Shintani \cite{shintania} introduced Dirichlet series whose coefficients
are the class numbers and studied extensively
as an example of {\em zeta functions of prehomogeneous vector spaces}.
This was generalized to over a general number field
using adelic language in \cite{wright}, and used to investigate relative
cubic extensions over the base field (\cite{dawra, dawrb})
or its integer ring (\cite{ddca}.)
Besides the rightmost pole at $s=1$,
the Dirichlet series has a mysterious second pole at $s=5/6$.
Our purpose is to study how those residues
behave when local conditions at finite places are imposed.

We prove theorems above in the following process.
In Section \ref{sec:notation} we review the notation
and invariant measure of \cite{ddca} those we use in this paper.
In Section \ref{sec:representation}
we recall and refine the parameterizations of cubic algebras 
by means of the space of binary cubic form
that we established in the previous paper \cite{ddca}.
In Section \ref{sec:zeta}
we introduce the global zeta function for the space of binary cubic forms.
We also introduce two partial zeta integrals,
those the contributions from irreducible forms and
from reducible forms. 
We use the result of Section \ref{sec:representation}
to express those integrals by Dirichlet series counting cubic algebras.

In Section \ref{sec:residue} we give the analytic continuations
and residue formulae of some Dirichlet series
including $\vartheta_{L_T}^\ird(\mathfrak a,s)$.
This give a proof of Theorem \ref{thm:mainthm_total} (1)
and Theorem \ref{thm:mainthm_stnz} (1)
with the values $\gA_{L_T}$, $\gB_{L_T}$
in Theorem \ref{thm:residueformula}.
The meromorphic continuations and residues
of the zeta integrals are obtained in \cite[Section 8]{ddca}.
The residues are expressed by means of local distributions and
hence the computation of the residues are reduced to the local theory.
The local theory for the space of binary cubic forms
were studied by Datskovsky-Wright \cite{dawra} in detail,
and their results make us the computation simple.
We note that the explicit formula \cite[Theorem 3.1]{dawra}
of the non-archimedean local zeta function plays an important role.

In Section \ref{sec:density} we prove
Theorem \ref{thm:mainthm_total} (2)
and Theorem \ref{thm:mainthm_stnz} (2).
Our tool to find density theorems is a
modified version \cite[Theorem 3]{sash}
of Landau's Tauberian theorem \cite[Hauptsatz]{landaub},
using a functional equation
to derive some informations on the error term.
To separate the contributions of reducible forms,
we also use the functional equations
of zeta functions for the space of binary quadratic forms.

\section{Notation and Invariant measures}
\label{sec:notation}

For notation we basically follow \cite{ddca}. But we have one exception.
If $V$ is a scheme define over a ring $R$ and $S$ is an $R$-algebra,
then in this paper
we denote the set of $S$-rational points of $V$ by $V(S)$,
not by $V_S$ as we denoted in \cite{ddca}.

For a finite set $X$ we denote by ${}^\# X$ its cardinality.
If an abstract group $G$ acts on a set $X$,
then for $x\in X$ we set ${\rm Stab}(G;x)=\{g\in G\mid gx=x\}$.
If $\mathfrak x\in G\backslash X$ is the class of $x\in X$,
we also denote the group by ${\rm Stab}(G;\mathfrak x)$,
which is well defined up to isomorphism.
The one-dimensional affine space is denoted by $\aff$.

Throughout this paper we fix a number field $k$.
We use the notation $\co$, $n$, $\zeta_k(s)$, $\Delta_k$,
$\gM_\infty$, $\gM_\R$, $\gM_\C$, $\gM_\fin$, $r_1$, $r_2$,
$k_v$, $q_v$, $\cl(k)$, $h_k$, $h_k^{(3)}$ and $\tau(\aaa)$
 as in Section \ref{sec:introduction}.
Let $\gM$ be the set of all places of $k$.
We put $\gc_k=\res_{s=1}\zeta_k(s)$.
For a fractional ideal $I$ of $k$, we denote by $N(I)$ its ideal norm.
The rings of adeles and finite adeles
are denoted by $\A$ and $\fa$.
We put $k_\infty=k\otimes_\Q\R$ and $\widehat\co=\co\otimes_\Z\widehat\Z$.
Note that $\widehat\co=\prod_{v\in\gMf}\co_v$,
$\fa=\widehat\co\otimes_\co k$,
and $\A=k_\infty\times \fa$.
The adelic absolute value $|\;|_\A$ on $\ma$ is normalized so that,
for $t\in\ma$, $|t|_\A$ is the module of multiplication by $t$
with respect to any Haar measure $dx$ on $\A$, i.e. $|t|_\A=d(tx)/dx$.
We define $|\;|_\infty$, $|\;|_{\fa}$ and $|\;|_v$
on $k_\infty^\times$, $\mfa$ and $k_v^\times$ similarly.
For a vector space $V$,
Let $\sS(V(\A))$, $\sS(V(k_v))$, $\sS(V(k_\infty))$ and $\sS(V(\fa))$
be the spaces of \bs s on each of the indicated domains.

For a fractional ideal $\aaa$, let $i(\aaa)\in\mfa(\subset \ma)$
be the corresponding idele, which is well defined up to
$\widehat\co^\times$-multiple.
That is, $i(\aaa)\in\mfa$ is characterized by
the condition $\aaa=k\cap i(\aaa)\widehat\co$.
Then $|i(\aaa)|_\A=N(\aaa)^{-1}$.
Notice that the infinite component of $i(\aaa)$ is trivial.
If there is no confusion we simply write $\aaa$ instead of $i(\aaa)$.
The set of characters of $\cl(k)$ is denoted by $\cl(k)^\ast$.
We regard $\omega\in\cl(k)^\ast$ as a character on
$\ma/\mk$ via the standard composition of the maps
$\ma/\mk\rightarrow
\ma/k_\infty^\times\mk\widehat\co^\times
	\cong \cl(k)\rightarrow \C^\times$.
Then $\omega(\aaa)=\omega(i(\aaa))$.

We give normalizations of invariant measures.
In general, for an algebraic group $X$ over $k$
with local measures $dx_v$ on $X(k_v)$ for all $v\in\gM$ are given,
then we always denote by
$dx_\infty	=\prod_{v\in\gM_\infty}dx_v$,
$dx_\fin	=\prod_{v\in\gMf}dx_v$ and
$d_\pr x	=\prod_{v\in\gM}dx_v$,
which are measures on $X(k_\infty)$, $X(\fa)$ and $X(\A)$, respectively.

For any  $v\in \gM_{\text{f}}$,
we choose a Haar measure $dx_v$ on $k_v$ 
to satisfy $\int_{\co_v}dx_v=1$.  
We write $dx_v$ for the ordinary Lebesgue measure if $v$ is real, 
and for twice the Lebesgue measure if $v$ is imaginary.  
For any $v\in \gM_{\text{f}}$, we normalize the Haar measure $\md t_v$
on $\mk_v$ such that $\int_{\co_v^{\times}} \md t_v = 1$.  
Let $\md t_v(x)=|x|_v^{-1}d x_v$ if $v\in\gM_\infty$.

Let $G=\gl_2$.
We review the normalization of the measure on $G(k_v)=\gl_2(k_v)$.
Let $B\subset G$ be the Borel subgroup
consisting of lower triangular matrices.
We normalize the right invariant measure
$db_v$ on $B(k_v)$ by
\[
\int_{B_{k_v}} f(b_v)db_v
=\int_{(k_v^\times)^2\times k_v}
	f\left(\twtw {s_v}00{t_v} \twtw 10{u_v}1\right)
	\left|\frac{t_v}{s_v}\right|_v \md s_v \md t_v du_v.
\]
Let $\cK=\prod_{v\in\gM}\cK_v$ where
$\cK_v={\rm O}(2),{\rm U}(2),\gl_2(\co_v)$ for
$v\in\gM_\R,\gM_\C,\gMf$, respectively.
We choose an invariant measure $d\kappa_v$ on $\cK_v$
such that $\int_{\cK_v}d\kappa_v=1$.
The group $G(k_v)$ has the decomposition
$G(k_v)=\cK_v B(k_v)$.
We choose an invariant measure on $G(k_v)$ by
$dg_v=d\kappa_v db_v$ for $g_v=\kappa_v b_v$.

Finally we express diagonal elements of $G$ as
$\diag(t_1,t_2)=\stwtw{t_1}00{t_2}$.

\section{The space of binary cubic forms and parameterization}
\label{sec:representation}
Let $G$ be the general linear group of rank $2$
and $V$ the space of binary cubic forms;
\begin{align*}
G&=\gl_2,\\
V&=
\{x=x(v_1,v_2)=x_0v_1^3+x_1v_1^2v_2+x_2v_1v_2^2+x_3v_2^3\mid x_i\in\aff\}.
\end{align*}
We define the action of $G$ on $V$ by
\[
(gx)(v)={(\det g)}^{-1}x(vg).
\]
The twist by $\det(g)^{-1}$ is to make the representation faithful.
For $x\in V$, let $P(x)$ be the discriminant;
\begin{equation*}
P(x)=x_1^2x_2^2-4x_0x_2^3-4x_1^3x_3+18x_0x_1x_2x_3-27x_0^2x_3^2.
\end{equation*}
Then we have $P(gx)=(\det g)^2P(x)$. We put $V^\sst=\{x\in V\mid P(x)\neq0\}$.

We first recall the parameterization of cubic algebras of $\co$.
Let $\aaa$ be a fractional ideal of $k$.
We put
\begin{align*}
G(k)\supset \gaa&=
\left\{\twtw abcd\ \vrule \ a\in\co,b\in\aaa,c\in\aaa^{-1},
d\in\co,ad-bc\in\co^\times\right\},\\
V(k)\supset \vaa&=
\{x\mid x_0\in \aaa, x_1\in \mathcal O,
x_2\in\aaa^{-1},x_3\in\aaa^{-2}\}.
\end{align*}
Then $\vaa$ is a $\gaa$-invariant submodule.
As in Section \ref{sec:introduction}, let
$\sC(\co)$ be the set of isomorphism classes 
of cubic algebras of $\co$.
For $R\in\sC(\co)$,
let $\aut(R)$ be the group of automorphisms of $R$ as an $\co$-algebra.
We put
$\sC(\co,\aaa)=\{R\in\sC(\co)\mid \st(R)=\aaa\}$.
In \cite[Section 3]{ddca} we establish the following
parameterization of $\sC(\co,\aaa)$.
\begin{prop}\label{prop:parameter_C(O,a)}
\begin{enumerate}[{\rm (1)}]
\item
There exists the canonical bijection
between $\sC(\co,\aaa)$ and $\gaa\backslash \vaa$
making the following diagram commutative:
\begin{equation*}
\xymatrix{
\gaa\backslash \vaa
\ar[r]^{\ \ \ \quad\quad\ \ \ } \ar[d]^{P}
& \sC(\co,\aaa) 
\ar[d]^{\text{discriminant}}\\
(\co^\times)^2\backslash \aaa^{-2}
\ar[r]^{\times\aaa^2\quad\quad}
&
\{\text{integral ideals of $\co$}\}.
}
\end{equation*}
Here, the right vertical arrow is to take the discriminant,
and the low horizontal arrow is given by multiplying $\aaa^2$.
For each $R\in\sC(\co,\aaa)$,
we denote by $x_R$ the corresponding element in $\gaa\backslash \vaa$
or its arbitrary representative in $\vaa$.
\item
We have $\aut(R)\cong \mathrm{Stab}(\gaa;x_R)$.
\end{enumerate}
\end{prop}

Let $K$ be either the number field $k$ or a local field $k_v$.
We next recall the geometric interpretation of the
$G(K)$-orbits in $V^\sst(K)$.
We denote by $\sC^{\sep}(K)$ the set of
isomorphism classes of separable cubic algebras of $K$.
For $x=x(v_1,v_2)\in V^\sst(K)$, we define
\[
Z_x	=\mathrm{Proj}\, K[v_1,v_2]/(x(v_1,v_2)),\qquad
K(x)	=\Gamma(Z_{x},\co_{Z_{x}}).
\]
We regard $K(x)$ as an element of $\sC^{\sep}(K)$.
It is well known that The map $x\mapsto K(x)$
gives a bijection between
$G(K)\backslash V^\sst(K)$
and $\sC^{\sep}(K)$.
By definition,
this correspondence is functorial with respect to the localization,
i.e., for $x\in V^\sst(k)$ and $v\in\gM$
we have $k_v(x)=k(x)\otimes k_v$.
We have the following compatibility
of the two parameterizations $x\mapsto x_R$ and $x\mapsto k(x)$.
\begin{lem}\label{lem:parameter_compatibility}
Let $R\in \sC(\co,\aaa)$ and $\Delta_{R/\co}\neq0$.
Then $x_R\in V^\sst(k)$ and $k(x_R)=R\otimes_\co k$.
\end{lem}
For $L_v\in \sC^{\sep}(k_v)$,
we put
\[
V_{L_v}=\{x\in V^\sst(k_v)\mid k_v(x)=L_v\},
\]
which is the $G(k_v)$-orbit in $V^\sst(k_v)$ corresponding to $L_v$.
This {\em should not be confused} to the set of $L_v$-rational points of $V$
nor the base change of $V$ to $L_v$.
Since we don't consider such objects in this paper
we hope this notation is not misleading.

\begin{defn}
Let $T$ be a set of places.
Let $L_T=(L_v)_{v\in T}$ be a $T$-tuple
where $L_v\in \sC^{\sep}(k_v)$ for each $v\in T$.
For such $L_T$, we define
\begin{align*}
V_{k,L_T}
&=	\{x\in V^\sst(k)\mid
		\text{$k(x)\otimes_k k_v=L_v$ for all $v\in T$}\},\\
\sC(\co,\aaa,L_T)
&=	\{R\in\sC(\co)\mid \text{$\st(R)=\aaa$,
		$R\otimes_\co k_v=L_v$ for all $v\in T$}\}.
\end{align*}
\end{defn}

\begin{prop}\label{prop:parameter_C(O,a,L_T)}
The correspondence in Proposition \ref{prop:parameter_C(O,a)}
induces a bijection between
$G_\aaa\backslash(V_\aaa\cap V_{k,L_T})$
and $\sC(\co,\aaa,L_T)$.
\end{prop}
\begin{proof}
Let $R\in\sC(\co,\aaa)$. By Lemma \ref{lem:parameter_compatibility},
$R\otimes_\co k_v=k(x_R)\otimes_k k_v$. Hence
$R\otimes_\co k_v=L_v$ if and only if $x_R\in V_{L_v}$.
Since $V_{k,L_T}=\bigcap_{v\in T}(V^\sst(k)\cap V_{L_v})$,
we have the bijection.
\end{proof}

Let $V^\sst(k)^\ird$ (resp.\ $V^\sst(k)^\rd$) be
the subset of $V^\sst(k)$ consisting of binary cubic forms
irreducible (resp.\ reducible) over $k$.
By the definition of $k(x)$,
$V^\sst(k)^\ird=\{x\in V^\sst(k)
	\mid\text{$k(x)/k$ is a cubic field extension}\}$.
Also we put
\begin{align*}
\sC(\co,\aaa,L_T)^\ird&=\{R\in\sC(\co,\aaa,L_T)\mid
	\text{$R\otimes_\co k$ is a cubic field extension of $k$}\},\\
\sC(\co,\aaa,L_T)^\rd&=\{R\in\sC(\co,\aaa,L_T)\mid
	\text{$R\otimes_\co k$ is a separable algebra and not a field}\}.
\end{align*}
Then by Lemma \ref{lem:parameter_compatibility},
we have the following refined version of
Proposition \ref{prop:parameter_C(O,a,L_T)}.
\begin{prop}\label{prop:parameter_C(O,a,L_T)_irdrd}
The correspondence of Proposition \ref{prop:parameter_C(O,a)}
induces a bijection between the following sets;
\begin{enumerate}[{\rm (1)}]
\item
$G_\aaa\backslash(V_\aaa\cap V_{k,L_T}\cap V^\sst(k)^\ird)$
and $\sC(\co,\aaa,L_T)^\ird$.
\item
$G_\aaa\backslash(V_\aaa\cap V_{k,L_T}\cap V^\sst(k)^\rd)$
and $\sC(\co,\aaa,L_T)^\rd$.
\end{enumerate}
\end{prop}

\section{Global zeta functions and Dirichlet series}\label{sec:zeta}
In this section we introduce the global zeta function
and related zeta integrals.
Using Propositions \ref{prop:parameter_C(O,a,L_T)}
and \ref{prop:parameter_C(O,a,L_T)_irdrd},
we see that these are integral expressions
of the Dirichlet series we want to investigate.

\begin{defn}
For $\Phi\in\sS(V_\A)$, $s\in\C$ and $\omega\in\cl(k)^\ast$,
we define the {\em global zeta function} by
\begin{equation*}
Z(\Phi,s,\omega)=\int_{G(\A)/G(k)}
|\det g|^{2s}\omega(\det g)
\sum_{x\in V^\sst(k)}\Phi(gx)d_\pr g.
\end{equation*}
Also for $\aaa\in\cl(k)$, we define
\begin{equation*}
Z_\aaa(\Phi,s)=\int_{G(k_\infty)G(\widehat\co)\diag(1,\aaa)G(k)/G(k)}
|\det g|^{2s}\sum_{x\in V^\sst(k)}\Phi(gx)d_\pr g.
\end{equation*}
\end{defn}
Recall the double coset decomposition
$G(\A)=
\coprod_{\aaa\in \cl(k)}G(k_\infty)G(\widehat\co)\diag(1,\aaa)G(k)$.
Since
$\omega(\det(G(k_\infty)G(\widehat\co)G(k)))
=\omega(k_\infty^\times\widehat\co^\times\mk)=1$,
we have the following.
\begin{lem}\label{lem:zetadecompose_idealclass}
\[
Z(\Phi,s,\omega)=\sum_{\aaa\in\cl(k)}\omega(\aaa)Z_\aaa(\Phi,s).
\]
\end{lem}
Let $T\supset \gM_\infty$ be a finite set of places.
For the rest of this paper we fix $T$ and $L_v\in \sC^\sep(k_v)$
for each $v\in T$. We put $L_T=(L_v)_{v\in T}$.
For our purpose, we also assume
$\Phi\in\sS(V(\A))$ is chosen as follows for the rest of this section.
\begin{asmp}\label{asmp:Phi}
We assume $\Phi$ is the product $\Phi=\prod_{v\in\gM}\Phi_v$,
where
\[
\text{$\Phi_v$ is}
\begin{cases}
\text{an arbitrary $\cK_v$-invariant function supported in $V_{L_v}$}
	&	v\in\gM_\infty,\\
\text{the characteristic function of $V(\co_v)\cap V_{L_v}$}
	&	v\in\gMf\cap T,\\
\text{the characteristic function of $V(\co_v)$}
	&	v\in\gMf\setminus T.
\end{cases}
\]
\end{asmp}
We put
$\Phi_\infty=\prod_{v\in\gM_\infty}\Phi_v$ and
$\Phi_\fin=\prod_{v\in\gMf}\Phi_v$.
Note that each $\Phi_v$ is $G(\co_v)$-invariant for $v\in\gMf$
and hence $\Phi_\fin$ is $G(\widehat\co)$-invariant.
\begin{defn}
We define
\begin{align*}
{\mathcal Z}_{L_T}(\Phi_\infty,s)
&=	\int_{G(k_\infty)}|P(g_\infty x)|_\infty^s
	\Phi_\infty(g_\infty x)dg_\infty\quad (x\in V_{k,L_T}),\\
\vartheta_{L_T}(\aaa,s)
&=	\sum_{x\in G_\aaa\backslash (V_\aaa\cap V_{k,L_T})}
\frac{({}^\#(\mathrm{Stab}(\gaa;x)))^{-1}}{N(\aaa)^{2s}|P(x)|_\infty^s}
=\sum_{R\in\sC(\co,\aaa,L_T)}
\frac{({}^\#\aut(R))^{-1}}{N(\Delta_{R/\co})^s}.
\end{align*}
\end{defn}
The function ${\mathcal Z}_{L_T}(\Phi_\infty,s)$ is called the
{\em local zeta function}.
Since $V_{k,L_T}\subset V(k_\infty)$ is contained in a
single $G(k_\infty)$-orbit $\prod_{v\in\gM_\infty}V_{L_v}$,
this does not depend on the choice of $x$.
Note that the second equality in the lower formula follows from
Propositions \ref{prop:parameter_C(O,a)} and
\ref{prop:parameter_C(O,a,L_T)}.
\begin{lem}\label{lem:unfolding_Z_a}
\[
Z_\aaa(\Phi,s)={\mathcal Z}_{L_T}(\Phi_\infty,s)\vartheta_{L_T}(\aaa,s).
\]
\end{lem}
\begin{proof}
Let $G(\widehat\co)_\aaa=\diag(1,\aaa)^{-1}G(\widehat\co)\diag(1,\aaa)$
and $\Phi_\aaa(x)=\Phi(\diag(1,\aaa)x)$.
Since the infinite part of $\aaa\in\ma$ is trivial,
the infinite part of $\Phi_\aaa$ coincides with $\Phi_\infty$.
We denote by $\Phi_{\fin,\aaa}$ the finite part of $\Phi_\aaa$,
hence $\Phi_\aaa=\Phi_\infty\times\Phi_{\fin,\aaa}$.
Then $\Phi_{\fin,\aaa}$ is $G(\widehat\co)_\aaa$-invariant.
Hence we have
\begin{align*}
Z_\aaa(\Phi,s)
&=\frac{1}{N(\aaa)^{2s}}
\int_{G(k_\infty) G(\widehat\co)_\aaa G(k)/G(k)}|\det g|_\A^{2s}
\sum_{x\in V^\sst(k)}\Phi_\aaa(gx)dg\\
&=\frac{1}{N(\aaa)^{2s}}
\int_{G(k_\infty) G(\widehat\co)_\aaa/G(k)\cap G(k_\infty) G(\widehat\co)_\aaa}
	|\det g_\infty|_\infty^{2s}
	\sum_{x\in V^\sst(k)}\Phi_{\fin,\aaa}(x)\Phi_\infty(g_\infty x)
	dg_\infty dg_\fin.
\end{align*}
For $x\in V^\sst(k)$, $\Phi_{\fin,\aaa}(x)=1$ if
\[
\diag(1,\aaa)x\in V(\widehat\co)
\qquad\text{and}\qquad
\diag(1,\aaa)_v x\in V_{L_v}\ \ \text{for all $v\in T$},
\]
and $0$ otherwise.
Here $\diag(1,\aaa)_v$ denote the $v$-component of $\diag(1,\aaa)$.
Since $V_{L_v}$ is a $G(k_v)$-orbit, the second condition holds
precisely when $x\in V_{L_v}$ for all $v\in T$.
Also we see $\diag(1,\aaa)^{-1}V(\widehat\co)\cap V(k)=V_\aaa$.
Hence $\Phi_{\fin,\aaa}(x)=1$ if $x\in V_\aaa\cap V_{k,L_T}$
and $0$ otherwise. Since
$G(k)\cap G(k_\infty) G(\widehat\co)_\aaa=G_\aaa$, we have
\[
Z_\aaa(\Phi,s)
=\frac{1}{N(\aaa)^{2s}}
\int_{G(k_\infty)/G_\aaa}
	|\det g_\infty|_\infty^{2s}
	\sum_{x\in V_\aaa\cap V_{L_T}}\Phi_\infty(g_\infty x)
	dg_\infty
\cdot
\int_{G(\widehat\co)_\aaa}dg_\fin.
\]
Since $G(\fa)$ is unimodular,
$\int_{G(\widehat\co)_\aaa}dg_\fin=1$.
Now the formula is obtained by the usual unfolding method.
Note that $|\det g_\infty|_\infty^2=|P(g_\infty x)/P(x)|_\infty$.
\end{proof}
Combined with Lemma \ref{lem:zetadecompose_idealclass}, we have
the following.
\begin{prop}\label{prop:unfolding_Z}
\[
Z(\Phi,s,\omega)=
	{\mathcal Z}_{L_T}(\Phi_\infty,s)
	\sum_{\aaa\in\cl(k)}\omega(\aaa)\vartheta_{L_T}(\aaa,s).
\]
\end{prop}
We next discuss the contributions of irreducible and reducible forms.
\begin{defn}
We put
\begin{align*}
Z^\ird(\Phi,s,\omega)
&=\int_{G(\A)/G(k)}
|\det g|^{2s}\tilde\omega(\det g)
\sum_{x\in V^\sst(k)^\ird}\Phi(gx)d_\pr g,\\
\vartheta_{L_T}^\ird(\aaa,s)
&
=\sum_{R\in\sC(\co,\aaa,L_T)^\ird}
\frac{({}^\#\aut(R))^{-1}}{N(\Delta_{R/\co})^s}.
\end{align*}
We define $Z^\rd(\Phi,s,\omega)$ and $\vartheta_{L_T}^\rd(\aaa,s)$ similarly.
\end{defn}
By Proposition \ref{prop:parameter_C(O,a,L_T)_irdrd},
the same argument of the proof of Lemma \ref{lem:unfolding_Z_a}
shows the following.
\begin{prop}\label{prop:unfolding_Z*}
\begin{align*}
Z^\ird(\Phi,s,\omega)
&=	{\mathcal Z}_{L_T}(\Phi_\infty,s)
	\sum_{\aaa\in\cl(k)}\omega(\aaa)\vartheta_{L_T}^\ird(\aaa,s),\\
Z^\rd(\Phi,s,\omega)
&=	{\mathcal Z}_{L_T}(\Phi_\infty,s)
	\sum_{\aaa\in\cl(k)}\omega(\aaa)\vartheta_{L_T}^\rd(\aaa,s).
\end{align*}
\end{prop}

\section{Analytic continuations and Residue formulae}\label{sec:residue}
In this section
we give proofs of Theorem \ref{thm:mainthm_total} (1)
and Theorem \ref{thm:mainthm_stnz} (1)
with the values $\gA_{L_T}$, $\gB_{L_T}$ in Theorem \ref{thm:residueformula}.
We also give the residue of $\vartheta_{L_T}^\rd(\aaa,s)$ at $s=1$.

We recall the residues of
$Z^\ird(\Phi,s,\omega)$ and $Z^\rd(\Phi,s,\omega)$.
For a while $\Phi\in\sS(V(\A))$ and $\Phi_v\in\sS(V(k_v))$ is arbitrary.
We introduce the following distributions.
\begin{defn}\label{defn:distribution_global}
For $\Phi\in\sS(V(\A))$, we define
\begin{align*}
\esA(\Phi)&=\int_{\A^4}\Phi(x_0,x_1,x_2,x_3)
		d_\pr x_0d_\pr x_1d_\pr x_2d_\pr x_3,\\
\esB(\Phi)&=\int_{\ma\times\A^3}
		|t|_\A^{1/3}\Phi(t,x_1,x_2,x_3)
		\md_\pr t d_\pr x_1 d_\pr x_2 d_\pr x_3,\\
\esC(\Phi)&=\int_{\ma\times\A^2}|t|_\A^2\Phi(0,t,x_2,x_3)
		\md_\pr t d_\pr x_2d_\pr x_3.
\end{align*}
\end{defn}
We have to give a notice for the definition of $\esB(\Phi)$,
because the integral defining $\esB(\Phi)$ itself does not converge.
(The integrals defining $\esA(\Phi)$ and $\esC(\Phi)$ converge.)
From the Iwasawa-Tate theory,
as a complex function of $s$ the integral
\[
\esB(\Phi,s)=\int_{\ma\times\A^3}
		|t|_\A^{s}\Phi(t,x_1,x_2,x_3)
		\md_\pr td_\pr x_1d_\pr x_2d_\pr x_3
\]
has meromorphic continuation to the whole complex plane
and is holomorphic except for $s=0,1$. The above definition means,
more precisely speaking, we put $\esB(\Phi)=\esB(\Phi,1/3)$.

The following is proved in \cite[Section 8]{ddca}.
\begin{thm}\label{thm:residueZ}
Let $\Phi$ be $\cK$-invariant.
The integrals $Z^\ird(\Phi,s,\omega)$ and $Z^\rd(\Phi,s,\omega)$
have meromorphic continuations to the whole complex plane.
Functions $(s-1)(s-5/6)Z^\ird(\Phi,s,\omega)$
and $(s-1)Z^\rd(\Phi,s,\omega)$
are holomorphic for $\re(s)>1/2$, and the residues of
$Z^\ird(\Phi,s,\omega)$ and $Z^\rd(\Phi,s,\omega)$
in this domain are given by the following;
\begin{align*}
\res_{s=1}Z^\ird(\Phi,s,\omega)
&=	\delta(\omega)
		2^{-1}\pi^{-r_1}(2\pi)^{-r_2}\gc_k \zeta_k(2)\esA(\Phi),\\
\res_{s=5/6}Z^\ird(\Phi,s,\omega)
&=	\delta(\omega^3)
		6^{-1}\Delta_k^{-1/2}\gc_k \esB(\Phi),\\
\res_{s=1}Z^\rd(\Phi,s,\omega)
&=	\delta(\omega)2^{-1}\gc_k \esC(\Phi).
\end{align*}
\end{thm}

We describe residues of $\vartheta_{L_T}^\ird(\aaa,s)$
and $\vartheta_{L_T}^\rd(\aaa,s)$
using this theorem and Proposition \ref{prop:unfolding_Z*}.
\begin{defn}\label{defn:distribution_local}
For $\Phi_v\in\sS(V(k_v))$, we define
\begin{align*}
\esA_v(\Phi_v)&=\int_{k_v^4}\Phi_v(x_0,x_1,x_2,x_3)dx_0dx_1dx_2dx_3,\\
\esB_v(\Phi_v)&=\int_{\mk_v\times k_v^3}
		|t|_v^{1/3}\Phi_v(t,x_1,x_2,x_3)\md tdx_1dx_2dx_3,\\
\esC_v(\Phi_v)&=\int_{\mk_v\times k_v^2}|t|_v^2\Phi_v(0,t,x_2,x_3)\md tdx_2dx_3.\end{align*}
Also we put
\[
\cZ_{L_v}(\Phi_v,s)
=	\int_{G(k_v)}|P(g_v x)|_v^s
	\Phi_v(g_v x)dg_v\quad (x\in V_{L_v}).
\]
\end{defn}

We now $\Phi=\prod_v\Phi_v$ is of the form in
Assumption \ref{asmp:Phi}.
\begin{defn}
Let $\Phi_v\in\sS(V(k_v))$ be as in Assumption \ref{asmp:Phi}.
\begin{enumerate}[{\rm (1)}]
\item
For $v\in\gM_\infty$, we define
\[
\alpha_v(L_v)=	\frac{\esA_v(\Phi_v)}{\cZ_{L_v}(\Phi_v,1)},
\quad
\beta_v(L_v)=	\frac{\esB_v(\Phi_v)}{\cZ_{L_v}(\Phi_v,5/6)},
\quad
\gamma_v(L_v)=	\frac{\esC_v(\Phi_v)}{\cZ_{L_v}(\Phi_v,1)}
\]
where $\Phi_v$ is chosen such that denominators do not vanish.
It is known that these do not depend on the choice of such $\Phi_v$.
\item
For $v\in\gMf\cap T$, let $\Phi_{v,0}$ be the characteristic
function of $V(\co_v)$. We put
\[
\alpha_v(L_v)=	\frac{\esA_v(\Phi_v)}{\esA_v(\Phi_{v,0})},
\quad
\beta_v(L_v)=	\frac{\esB_v(\Phi_v)}{\esB_v(\Phi_{v,0})},
\quad
\gamma_v(L_v)=	\frac{\esC_v(\Phi_v)}{\esC_v(\Phi_{v,0})}.
\]
\item
Further we let the product as follows;
\[
\alpha_T(L_T)=\prod_{v\in T}\alpha_v(L_v),
\quad
\beta_T(L_T)=\prod_{v\in T}\beta_v(L_v),
\quad
\gamma_T(L_T)=\prod_{v\in T}\gamma_v(L_v).
\]
\end{enumerate}\end{defn}
\begin{thm}\label{thm:dirichlet}
\begin{enumerate}[{\rm (1)}]
\item
The Dirichlet series $\vartheta_{L_T}^\ird(\aaa,s)$
and $\vartheta_{L_T}^\rd(\aaa,s)$
have meromorphic continuations to the whole complex plane.
Functions $(s-1)(s-5/6)\vartheta_{L_T}^\ird(\aaa,s)$
and $(s-1)\vartheta_{L_T}^\rd(\aaa,s)$
are holomorphic for $\re(s)>1/2$.
\item
The residues of $\vartheta_{L_T}^\ird(\aaa,s)$
and $\vartheta_{L_T}^\rd(\aaa,s)$
in the region $\re(s)>1/2$
are given by;
\begin{align*}
\res_{s=1}\vartheta_{L_T}^\ird(\aaa,s)
&=	h_k^{-1}2^{-1}\pi^{-r_1}(2\pi)^{-r_2}\gc_k \zeta_k(2)\alpha_T(L_T),\\
\res_{s=5/6}\vartheta_{L_T}^\ird(\aaa,s)
&=	\tau(\aaa)h_k^{(3)}h_k^{-1}
		6^{-1}\Delta_k^{-1/2}\gc_k\zeta_k(1/3)\beta_T(L_T),\\
\res_{s=1}\vartheta_{L_T}^\rd(\aaa,s)
&=	h_k^{-1}2^{-1}\gc_k\zeta_k(2)\gamma_T(L_T).
\end{align*}
\end{enumerate}
\end{thm}
\begin{proof}
By the theory of prehomogeneous vector spaces,
the local zeta function ${\mathcal Z}_{L_T}(\Phi_\infty,s)$
has meromorphic continuation to the whole complex plane.
Moreover, for any $s\in\C$ such that $\re(s)>1/6$,
we can choose $\Phi_\infty$ such that ${\mathcal Z}_{L_T}(\Phi_\infty,s)\neq0$.
Hence (1) follows from Proposition \ref{prop:unfolding_Z*},
Theorem \ref{thm:residueZ} and the orthogonality of characters.

We consider (2). We prove the second formula.
The first and third formulae are proved similarly.
Let
$\vartheta_{L_T}^\ird(\aaa,s,\omega)=
\sum_{\aaa\in\cl(k)}\omega(\aaa)\vartheta_{L_T}^\ird(\aaa,s)$,
which is the Dirichlet series appeared in 
the first formula of Proposition \ref{prop:unfolding_Z*}.
For $\Psi_\fin\in\sS(V(\fa))$ we put the ``finite part''
$\esB_\fin(\Psi_\fin)$ of $\esB$ by
evaluating at $s=1/3$ of the analytic function
$
\int_{\mfa\times\fa^3}
		|t|_{\fa}^{s}\Psi_\fin(t,x_1,x_2,x_3)
			\md_\fin td_\fin x_1d_\fin x_2d_\fin x_3$.
Then $\esB(\Phi)=\esB_\fin(\Phi_\fin)\prod_{v\in\gM_\infty}\esB_v(\Phi_v)$.
Let $\Phi_{\fin,0}$ be the characteristic function of $V(\widehat\co)$.
Then we have $\esB(\Phi_{\fin,0})=\zeta_k(1/3)$ and
$\esB_\fin(\Phi_{\fin})/\esB_\fin(\Phi_{\fin,0})
=\prod_{v\in T\cap\gMf} \esB_v(\Phi_{v})/\esB_v(\Phi_{v,0})$
by the analytic continuations.
Hence by Proposition \ref{prop:unfolding_Z*} and Theorem \ref{thm:residueZ}
we have
\begin{align*}
\res_{s=5/6}\vartheta_{L_T}^\ird(\aaa,s,\omega)
&=\delta(\omega^3)\frac{\gc_k\zeta_k(1/3)}{6\Delta_k^{1/2}}
\frac{\prod_{v\in\gM_\infty}{\esB_v(\Phi_v)}}
{{\mathcal Z}_{L_T}(\Phi_\infty,5/6)}
\prod_{v\in T\cap \gMf}\frac{\esB_v(\Phi_v)}{\esB_v(\Phi_{v,0})}\\
&=\delta(\omega^3)6^{-1}\Delta_k^{-1/2}\gc_k\zeta_k(1/3)\beta_T(L_T).
\end{align*}
By the orthogonality of the characters,
$\vartheta_{L_T}^\ird(\aaa,s)=\sum_{\omega\in\cl(k)^\ast}
\omega(\aaa)^{-1}\vartheta_{L_T}^\ird(\aaa,s,\omega)$.
From the identity
$\sum_{\omega\in\cl(k)^\ast}
\omega(\aaa)^{-1}\delta(\omega^3)=\tau(\aaa)h_k^{(3)}/h_k$,
we have the second formula of (2).
\end{proof}
From now on we consider $\alpha_v(L_v)$, $\beta_v(L_v)$, and $\gamma_v(L_v)$.
For $v\in\gM_\infty$, this is done by
Shintani \cite{shintania} and Datskovsky-Wright \cite{dawra}
and we already used their result in \cite{ddca}.
For the convenience of the reader, we give these values
in Table \ref{table:infinite}:
\begin{prop}\label{prop:infinite_value}
For $v\in\gM_\infty$,
$\alpha_v(L_v)$, $\beta_v(L_v)$ and $\gamma_v(L_v)$
are given by Table \ref{table:infinite}.
\end{prop}
\begin{table}
\[
\begin{array}{|c|c||c|c|c|}
\hline
v	&L_v
	& \alpha_v(L_v)
	& \beta_v(L_v)
	& \gamma_v(L_v)\\[1mm]
\hline
\hline
\R	&\R\times\R\times\R
	&\dfrac{\pi}{6}
	&\dfrac{\sqrt3\Gamma(1/3)^3}{4\pi}
	&\dfrac12\\[3mm]
\hline
\R	&\R\times\C
	&\dfrac{\pi}{2}
	&\dfrac{3\Gamma(1/3)^3}{4\pi}
	&\dfrac12\\[3mm]
\hline
\C	&\C\times\C\times\C
	&\dfrac{2\pi}{6}
	&\dfrac{\sqrt3\Gamma(1/3)^6}{8\pi^2}
	&\dfrac12\\[3mm]
\hline
\end{array}
\]
\vspace*{5pt}
\caption{}%
\label{table:infinite}
\end{table}

We compute $\alpha_v(L_v)$, $\beta_v(L_v)$, and $\gamma_v(L_v)$ for $v\in\gMf$
using the local theory developed by Datskovsky-Wright \cite{dawra},
especially Theorems 3.1, 5.1, 5.2
and Propositions 5.1, 5.3.
We fix $v\in\gMf$ for the rest of this section.
To simplify the notation, we drop the subscript $v$
and write $L$, $q$, $\alpha$, $\beta$, $\gamma$, $\esA$, $\esB$ and $\esC$
instead of writing $L_v$, $q_v$, $\alpha_v$, $\beta_v$, $\gamma_v$
$\esA_v$, $\esB_v$ and $\esC_v$.
To stress the dependence $\Phi_v$ in Assumption \ref{asmp:Phi} on $L$,
we write $\Phi_v=\Phi_L$.
As in Section \ref{sec:introduction},
we denote by $\theta_L$
the order of the automorphisms of $L$ as a $k_v$-algebra.
Let $\zeta_L(s)$ be the local zeta function associated with
the $k_v$-algebra $L$. Namely, we define
\[
\zeta_L(s)=
\begin{cases}
(1-q^{-s})^{-3}
&	L=k_v\times k_v\times k_v,\\
(1-q^{-s})^{-1}(1-q^{-2s})^{-1}
&	L=k_v\times\text{(quad. unramified ext. of $k_v$)},\\
(1-q^{-3s})^{-1}
&	L=\text{(cubic unramified ext. of $k_v$)},\\
(1-q^{-s})^{-2}
&	L=k_v\times\text{(quad. ramified ext. of $k_v$)},\\
(1-q^{-s})^{-1}
&	L=\text{(cubic ramified ext. of $k_v$)}.\\
\end{cases}
\]
Let $\Delta_{L/k_v}$ be the discriminant of $L/k_v$
and $\Delta_L$ its norm.
An element $x_L$ in the $G(k_v)$-orbit $V_{L}$ is called a
{\em standard orbital representative} if $x_L\in V(\co_v)$ and
$P(x_L)\in\co_v$ generates the ideal $\Delta_{L/k_v}$.
It is easy to see that such an element exists for arbitrary $L$.
We fix such $x_L$.
We define
\[
\Omega_L(\Phi_L,s)=\int_{G(k_v)}|\det g_v|_v^{2s}\Phi_L(g_v x_L)dg_v,
\]
which equals to $|P(x_L)|_v^{-s}\cZ_{L_v}(\Phi_v,s)$
and hence do not depend on the choice of $x_L$.
This function plays an important role in the computation.

Before starting the computation,
we will compare and adjust the notation
in \cite{dawra} to ours.
In \cite{dawra} they denote by $A$ the set of $G(k_v)$-orbits of $V^\sst(k_v)$
and by $\alpha$ any of its element.
Hence there is the canonical correspondence between
their $A$ and our $\sC^\sep(k_v)$.
If $L\in \sC^\sep(k_v)$ corresponds to $\alpha\in A$, then
$\Omega_L(\Phi_L,s)$ equals to what they denoted
by $Z_\alpha(\omega_{2s},\Phi_L)$.
(They defined the distribution $Z_\alpha(\omega,\Phi)$ in p.39.)
Also the value $c_\alpha$ they introduced in p.38
equals to $\theta_L^{-1}\Delta_L^{-1}$, since
their $o(\alpha)$ equals to $\theta_L$ and
their $x_\alpha$ is also a standard orbital representative
for the orbit $V_{L}$.

The following beautiful formula
is a variation of \cite[Theorem 3.1]{dawra}.
\begin{lem}\label{lem:explicit_localzeta}
\[
\Omega_L(\Phi_L,s)=(1-q^{-4s})^{-1}(1-q^{-6s+1})^{-1}\zeta_L(2s)\zeta_L(4s)^{-1}.
\]
\end{lem}
\begin{proof}
Let $\Phi_0\in\sS(V(k_v))$ be the characteristic function of $V(\co_v)$.
Then since $\Phi_L$ is the characteristic function of
$V(\co_v)\cap G(k_v)x_L$, $\Phi_L(g_v x_L)=\Phi_0(g_v x_L)$
for all $g_v\in G(k_v)$. Hence
\[
\Omega_L(\Phi_L,s)=\int_{G(k_v)}|\det g_v|_v^{2s}\Phi_0(g_v x_L)dg_v.
\]
The right hand side is what Datskovsky and Wright
gave the explicit formula in \cite[Theorem 3.1]{dawra}.
Since their notation $\omega(\pi)$ is $q^{-2s}$ in our setting,
we have the formula.
\end{proof}

\begin{prop}\label{prop:alpha_v}
\[
\alpha(L)
=	\frac{(1-q^{-1})(1-q^{-2})}{(1-q^{-4})(1-q^{-5})}
	\cdot	\theta_L^{-1}\Delta_L^{-1}
	\cdot	\frac{\zeta_L(2)}{\zeta_L(4)}.
\]
\end{prop}
\begin{proof}
Obviously $\esA(\Phi_0)=1$ and hence $\alpha(L)=\esA(\Phi_L)$.
This is expressed by $q^{2e}\widehat{\Phi_L}(0)$ in the notation of
\cite{dawra}.
Note that their additive measure on $k_v$ is $q^{-e/2}$ times to ours.
Since $\Phi_L(x)=0$ if $x\not\in V_L$,
by \cite[Propostion 5.1]{dawra} we have
\[
\esA(\Phi_L)=(1-q^{-1})(1-q^{-2})\theta_L^{-1}\Delta_L^{-1}\Omega_L(\Phi_L,1).
\]
Hence from Lemma \ref{lem:explicit_localzeta} we have the formula.
\end{proof}

\begin{prop}\label{prop:beta_v}
\[
\beta(L)
=	\frac{(1-q^{-1/3})(1-q^{-1})}{(1-q^{-10/3})(1-q^{-4})}
	\cdot	\theta_L^{-1}\Delta_L^{-1}
	\cdot	\frac{\zeta_L(1/3)\zeta_L(5/3)}{\zeta_L(2/3)\zeta_L(10/3)}.
\]
\end{prop}
\begin{proof}
Since $\esB(\Phi_0)=(1-q^{-1/3})^{-1}$ we have
$\beta(L)=(1-q^{-1/3})\esB(\Phi_L)$.
If we use the notation of \cite{dawra},
$\esB(\Phi_L)$ equals to $q^{3e/2}\Sigma_4(1,\Phi_L)$
where $1$ denotes the trivial character on $k_v$.
Hence if $L$ corresponds to $\alpha\in A$ in their notation,
by \cite[Theorem 5.2]{dawra}
we have
\[
\esB(\Phi_L)=(1-q^{-1})\theta_L^{-1}\Delta_L^{-1}a_\alpha(1)\Omega_L(\Phi_L,5/6).
\]
The value $a_\alpha(1)$ is evaluated in \cite[Proposition 5.3]{dawra}
and is found to be $\zeta_L(1/3)\zeta_L(2/3)^{-1}$.
Now Lemma \ref{lem:explicit_localzeta}
gives the desired formula.
\end{proof}

Theorem \ref{thm:mainthm_total} (1)
and Theorem \ref{thm:mainthm_stnz} (1)
with the values $\gA_{L_T}$, $\gB_{L_T}$ in Theorem \ref{thm:residueformula}
follows from Theorem \ref{thm:dirichlet} and
Propositions \ref{prop:infinite_value},
\ref{prop:alpha_v}, \ref{prop:beta_v}.
We also evaluate $\gamma(L)$.
This gives the explicit formula of the residue of $\vartheta_{L_T}^\rd(\aaa,s)$
at $s=1$.

\begin{prop}\label{prop:gamma_v}
\[
\gamma(L)
=
\begin{cases}
\dfrac{(1-q^{-1})(1-q^{-2})}{(1-q^{-4})(1-q^{-5})}
	\cdot	2^{-1}\Delta_L^{-1}
	\cdot	\dfrac{\zeta_L(2)}{\zeta_L(4)}
&	\text{$L$ is not a field},\\
0
&	\text{$L$ is a field}.\\
\end{cases}
\]
\end{prop}
\begin{proof}
The value $\esC(\Phi_L)$ equals to $q^e\Sigma_2(\Phi_L)$
if we use the notation of \cite{dawra}.
Hence by \cite[Theorem 5.2]{dawra},
this equals to $(1-q^{-1})2^{-1}\Delta_L^{-1}\Omega_L(\Phi_L,1)$
if $L$ is not a field and $0$ otherwise.
Since $\esC(\Phi_0)=(1-q^{-2})^{-1}$, we have the formula.
\end{proof}

\begin{rem}\label{rem:sum_localdensity}
Let $\sE_i^{\rm\, rm}(k_v)$ the set of isomorphism
classes of totally ramified extensions of $k_v$ of degree $i$.
For $F\in \sE_i^{\rm\, rm}(k_v)$, let
$\Delta_F$ be the norm of the relative discriminant of $F/k_v$
and $\theta_F$ the order of automorphisms of $F$.
In \cite{serrec} Serre established a beautiful formula
\[
\sum_{F\in\sE_i^{\rm\, rm}(k_v)}
	\theta_F^{-1}\Delta_F^{-1}
=	q^{-i+1}.
\]
On the other side, by gathering all the splitting type
at a fixed $v\in\gMf$, we have
\[
\sum_{L\in\sC^\sep(k_v)}\alpha(L)
=\sum_{L\in\sC^\sep(k_v)}\beta(L)
=\sum_{L\in\sC^\sep(k_v)}\gamma(L)
=1.
\]
By computation we see that this is equivalent to the Serre's
formula for $i=2$ and $3$.
\end{rem}

\section{Density theorems}\label{sec:density}
In this section we give proofs of
Theorem \ref{thm:mainthm_total} (2) and
Theorem \ref{thm:mainthm_stnz} (2).
Throughout this section we fix $L_T$.
Recall that we define Dirichlet series
$\vartheta_{L_T}(\aaa,s)$,
$\vartheta_{L_T}^\ird(\aaa,s)$,
$\vartheta_{L_T}^\rd(\aaa,s)$ in Section \ref{sec:zeta}.
These satisfy
$\vartheta_{L_T}(\aaa,s)
=\vartheta_{L_T}^\ird(\aaa,s)+
\vartheta_{L_T}^\rd(\aaa,s)$.
We put
$\vartheta_{L_T}(s)=\sum_{\aaa\in\cl(k)}\vartheta_{L_T}(\aaa,s)$.
Similarly we define $\vartheta_{L_T}^\ird(s)$, $\vartheta_{L_T}^\rd(s)$
as the sum of the Dirichlet series over $\aaa\in\cl(k)$.
We write
\[\textstyle
\vartheta_{L_T}(s)=\sum_{m\geq1}a_m/m^{s},\quad
\vartheta_{L_T}^\ird(s)=\sum_{m\geq1}a_m^\ird/m^s,\quad
\vartheta_{L_T}^\rd(s)=\sum_{m\geq1}a_m^\rd/m^s,
\]
so that $a_n=a_n^\ird+a_n^\rd$.
Until the proof of Lemma \ref{lem:a_mrd} we prove the following.
\begin{prop}\label{prop:a_mird}
For any $\varepsilon>0$,
\[
\sum_{m<X}a_m^\ird
	=\gA_{L_T}X+(5/6)^{-1}\gB_{L_T}X^{5/6}
		+O(X^{\frac{5n-1}{5n+1}+\varepsilon})
\qquad (X\to\infty).
\]
\end{prop}
To prove the proposition we first give an estimate of
the function $\sum_{m<X}a_m$.
We put $\gC_{L_T}=2^{-1}\gc_k\zeta_k(2)\prod_{v\in T}\gamma_v(L_v)$,
which is the residue of $\vartheta_{L_T}^\rd(s)$ at $s=1$. 
\begin{lem}\label{lem:a_m}
For any $\varepsilon>0$,
\[
\sum_{m<X}a_m
	=(\gA_{L_T}+\gC_{L_T})X+(5/6)^{-1}\gB_{L_T}X^{5/6}
		+O(X^{\frac{4n-1}{4n+1}+\varepsilon})
\qquad (X\to\infty).
\]
\end{lem}
\begin{proof}
In \cite{wright} Wright proved that
the global zeta function $Z(\Phi,s,\omega)$
can be continued holomorphically to the whole complex plane
except for possible simple poles at $s=0,1/6,5/6,1$ and
satisfies the functional equation
\[
Z(\Phi,s,\omega)=Z(\hat\Phi,1-s,\omega^{-1})
\]
where $\hat\Phi$ is a suitable Fourier transform of $\Phi$.
Let $\omega\in\cl(k)^\ast$ be the trivial character and
$\Phi\in\sS(V(\A))$ as in Assumption \ref{asmp:Phi}.
Then by Proposition \ref{prop:unfolding_Z}
$Z(\Phi,s,\omega)={\mathcal Z}_{L_T}(\Phi_\infty,s)\vartheta_{L_T}(s)$.
Since $\hat\Phi$ is also a $G(\widehat\co)$-invariant function,
we have a similar decomposition for $Z(\hat\Phi,s,\omega^{-1})$.
These combined with archimedean local theory show
that $\vartheta_{L_T}(s)$ is holomorphic except for $s=1,5/6$ and
satisfies a functional equation of the form
\[
\vartheta_{L_T}(1-s)=\Gamma(s)^{2n}\Gamma(s-\frac16)^n\Gamma(s+\frac16)^n
\sum_{\lambda\in\Lambda}p_\lambda(e^{\pi\sqrt{-1}s/2},e^{-\pi\sqrt{-1}s/2})
\xi_\lambda(s)
\]
where $\Lambda$ is a finite index set and for each $\lambda\in\Lambda$,
$p_\lambda(x,y)$ is a polynomial in $x,y$ of degree less than $4n$
and $\xi_\lambda(s)$ is a Dirichlet series with absolute convergence domain
$\re(s)>1$.
By Theorem \ref{thm:dirichlet} the residues of $\vartheta_{L_T}(s)$
at $s=1$ and $5/6$
are $\gA_{L_T}+\gC_{L_T}$ and $\gB_{L_T}$, respectively.
Hence the proposition follows from
the modified Landau theorem \cite[Theorem 3]{sash}.
\end{proof}
This lemma reduces the proof of Proposition \ref{prop:a_mird}
to an estimate of $\sum_{m<X}a_m^\rd$.
To give an estimate of this function we introduce another \pv.
Let $B=\br_2\subset G$ and $W$ the subspace of $V$ having a linear factor $v$;
\[
B=\left\{\twtw *0**\right\},\quad
W=\{y=y(u,v)=v(y_1u^2+y_2uv+y_3v^2)\mid y_1,y_2,y_3\in\aff\}.
\]
Then $W$ is an invariant subspace of $B$ and $(B,W)$
is also a prehomogeneous vector space.
As in \cite{ddca} , we use the 
functional equation of the zeta function for this space
to study $\sum_{m<X}a_m^\rd$.
(This idea is due to Shintani \cite{shintanib}.)
For $\aaa,\ccc\in\cl(k)$, we put
\begin{align*}
B(k)\supset \bac&=
\left\{\twtw t0up\ \vrule\ 
	t,p\in\co^\times,u\in\aaa^{-1}\ccc^{-2}\right\},\\
W(k)\supset \wac&=
\{y\mid y_1\in\ccc, y_2\in\aaa^{-1}\ccc^{-1},y_3\in\aaa^{-2}\ccc^{-3}\}.
\end{align*}
Then $\wac$ is $\bac$-invariant.
Let $W_{k,L_T}=W(k)\cap V_{k,L_T}$.
We define
\[
\eta_{L_T}(s)=\sum_{\aaa,\ccc\in\cl(k)}
\sum_{y\in \bac\backslash (\wac\cap W_{k,L_T})}
\frac{({}^\#{\rm Stab}(\bac;y))^{-1}}{N(\aaa)^{2s}|P(y)|_\infty^s}.
\]
We write $\eta_{L_T}(s)=\sum_{m\geq1}b_m/m^s$.
Note that for $y\in\wac$, $P(y)\in\aaa^{-2}$ and hence
$N(\aaa)^2|P(y)|_\infty$ is an integer.
Then we have the following.
\begin{lem}\label{lem:b_m}
For any $\varepsilon>0$,
\[
\sum_{m<X}b_m
	=\gC_{L_T}X+O(X^{\frac{5n-1}{5n+1}+\varepsilon})
\qquad (X\to\infty).
\]
\end{lem}
\begin{proof}
We give an integral expression of $\eta_{L_T}(s)$.
For $\Psi\in\sS(W(\A))$ and $s\in\C$,
define
\[
Y(\Psi,s)=\int_{B(\A)/B(k)}|\det b|_\A^{2s}\sum_{y\in W(k)\cap V^\sst(k)}\Psi(by)db.
\]
Also for $\Psi_\infty\in\sS(W(k_\infty))$, we put
\[
\mathcal Y_{L_T}(\Psi_\infty,s)
=\int_{B(k_\infty)}|P(b_\infty y)|_\infty^s\Psi_\infty(b_\infty y)db_\infty
\quad (y\in W_{k,L_T}.)
\]
Let $\Phi\in\sS(V(\A))$ be as in Assumption \ref{asmp:Phi}.
We regard $\Phi\in\sS(W(\A))$ via the pullback of
the inclusion $W(\A)\rightarrow V(\A)$.
Then since
\begin{gather*}
B(\A)=\coprod_{\aaa,\ccc\in\cl(k)}
B(k_\infty)B(\widehat\co)\cdot \diag(\ccc,\aaa^{-1}\ccc^{-1})\cdot B(k),\\
\bac=
B(k)\cap \diag(\ccc,\aaa^{-1}\ccc^{-1})^{-1}\cdot B(k_\infty)B(\widehat\co)
\cdot\diag(\ccc,\aaa^{-1}\ccc^{-1}),\\
\wac=
W(k)\cap(\diag(\ccc,\aaa^{-1}\ccc^{-1})^{-1}W(\widehat\co)\times W(k_\infty)),
\end{gather*}
by the similar unfolding method as in Proposition \ref{prop:unfolding_Z}
we have
\[
Y(\Phi,s)=\mathcal Y_{L_T}(\Phi_\infty,s)\eta_{L_T}(s).
\]
We recall from \cite{ddca}
some results on zeta functions for $(B,W)$.
The residue of $Y(\Phi,s)$ at $s=1$ is exactly the same
as that of $Z^\rd(\Phi,s)$. Also since $\Phi_\infty$ is
$\prod_{v\in\gM_\infty}\cK_v$-invariant,
the local zeta functions
${\mathcal Z}_{L_T}(\Phi_\infty,s)$ and $\mathcal Y_{L_T}(\Phi_\infty,s)$
also coincide.
Hence the residue of $\eta_{L_T}(s)$ at $s=1$ is $\gC_{L_T}$.
Moreover, the Dirichlet series
$H(s)=\eta_{L_T}(s)\zeta_k(4s)$ is entire after multiplied by
$(s-1/2)^3(s-1)$, and it satisfies a functional equation of the form
\[
H(1-s)
=\Gamma(2s-1)^{2n}\Gamma(s)^n\Gamma(3s-3/2)^n
\sum_{\lambda\in\Lambda}q_\lambda(e^{\pi\sqrt{-1}s/2},e^{-\pi\sqrt{-1}s/2})
\theta_\lambda(s)
\]
where $\Lambda$ is a finite index set and for each $\lambda\in\Lambda$,
$q_\lambda(x,y)$ is a polynomial in $x,y$ of degree less than $8n$
and $\theta_\lambda(s)$ is a Dirichlet series with absolute convergence domain
$\re(s)>1$. Hence by the same argument as in \cite[Proposition 7.18]{ddca}
we have the lemma.
\end{proof}
We introduce the following notation.
\begin{defn}
Let $k(s)=\sum_{m\geq1}k_m/m^s$ and $l(s)=\sum_{m\geq1}l_m/m^s$
be Dirichlet series having absolute convergence domains.
We say $k(s)$ is bounded by $l(s)$ if $|k_m|\leq l_m$ for all $m\geq1$,
and write $k(s)\preceq l(s)$ in this situation.
\end{defn}
Using Lemma \ref{lem:b_m} we have the following estimate,
which proves Proposition \ref{prop:a_mird}.
Notice that for a positive sequence $\{a_m\}$ and
a positive constant $\rho$,
the series $\sum_{m\geq1}a_m/m^s$ converges for $\re(s)>\rho$
if and only if
$\sum_{m<X}a_n=O(X^{\rho+\epsilon})$
for any $\epsilon>0$.

\begin{lem}\label{lem:a_mrd}
For any $\varepsilon>0$,
\[
\sum_{m<X}a_m^\rd
	=\gC_{L_T}X+O(X^{\frac{5n-1}{5n+1}+\varepsilon})
\qquad (X\to\infty).
\]
\end{lem}
\begin{proof}
Let $V_k^{(1)}=\{x\in V^\sst(k)\mid k(x)=k\times k\times k\}$,
$V_k^{(2)}=V^\sst(k)\setminus (V^\sst(k)^\ird\cup V_k^{(1)})$
and put
$\vaa^{(i)}=\vaa\cap V_k^{(i)}$,
$\wac^{(i)}=\wac\cap V_k^{(i)}$ for $i=1,2$.
In \cite[Proposition 3.12]{ddca}, we construct a bijective map
\[
\psi_\aaa\colon\coprod_{\ccc\in\cl(k)}(\bac\backslash\wac^{(2)})
	\longrightarrow\gaa\backslash\vaa^{(2)}
\]
which preserves the value of $P$ up to $(\co^\times)^2$-multiple
and $\mathrm{Stab}(\bac;y)
\cong\mathrm{Stab}(\gaa;\psi_{\aaa}(y))$ for any $y\in \wac^{(2)}$.
From the construction it is easy to see that this induces a bijective map
\[
\psi_{\aaa,L_T}\colon
\coprod_{\ccc\in\cl(k)}
	\left(\bac\backslash(\wac^{(2)}\cap W_{k,L_T})\right)
	\longrightarrow\gaa\backslash(\vaa^{(2)}\cap V_{k,L_T}).
\]
On the other hand by the local theory, the Dirichlet series
\[
\sum_{x\in \gaa\backslash(\vaa^{(1)}\cap V_{k,L_T})}
			N(\aaa)^{-2s}|P(x)|_\infty^{-s},\quad
\sum_{y\in \bac\backslash(\wac^{(1)}\cap W_{k,L_T})}
			N(\aaa)^{-2s}|P(y)|_\infty^{-s}
\]
are bounded by the Dirichlet series
$\zeta_k(2s)^3\zeta_k(6s-1)/\zeta_k(4s)^2$.
Especially these Dirichlet series converges for $\re(s)>1/3$.
Hence $\sum_{m<X}a_m^\rd=\sum_{m<X}b_m+O(X^{1/3+\epsilon})$
where $\epsilon>0$ is arbitrary and this finishes the proof.
\end{proof}
We now give a proof of Theorem \ref{thm:mainthm_total} (2).
\begin{thm}\label{thm:h_L_T}
 Let $h_{L_T}(X)$ be as in Section \ref{sec:introduction}.
For any $\varepsilon>0$,
\[
h_{L_T}(X)
	=\gA_{L_T}X+(5/6)^{-1}\gB_{L_T}X^{5/6}
		+O(X^{\frac{5n-1}{5n+1}+\varepsilon})
\qquad (X\to\infty).
\]
\end{thm}
\begin{proof}
Let
\[
\tilde\vartheta_{L_T}^\ird(s)
=\sum_{n\geq1}\frac{\tilde a_n^\ird}{n^s}
=	\sum_{R\in\sC(\co,L_T)^\ird}
		\frac{1}{N(\Delta_{R/\co})^s}.
\]
Then $h_{L_T}(X)=\sum_{m<X}\tilde a_m^\ird$.
We compare $\vartheta_{L_T}^\ird(s)$ and $\tilde\vartheta_{L_T}^\ird(s)$.
Let $\gG$ be the set of isomorphism classes of cyclic cubic extensions of $k$.
Then $\aut(R)$ is non-trivial only if $R\otimes k\in\gG$. Hence
as in the proof of \cite[Theorem 7.20]{ddca},
\begin{align*}
\tilde\vartheta_{L_T}^\ird(s)-
\vartheta_{L_T}^\ird(s)
&\preceq\sum_{F\in \gG}\sum_{R\in\sC(\co),R\otimes k=F}
		\frac{1}{N(\Delta_{R/\co})^s}\\
&=	\sum_{F\in \gG}N(\Delta_{F/k})^{-s}
		\zeta_k(4s)\zeta_k(6s-1)\zeta_F(2s)\zeta_F(4s)^{-1}\\
&\preceq\zeta_k(2s)^3\zeta_k(6s-1)\zeta_k(4s)^{-2}
	\sum_{F\in \gG}N(\Delta_{F/k})^{-s}
\end{align*}
and the Dirichlet series in the
last term has absolute convergence domain $\re(s)>1/2$.
This shows that
$h_{L_T}(X)=\sum_{m<X}\tilde a_m^\ird=\sum_{m<X}a_m^\ird+O(X^{1/2+\epsilon})$
where $\epsilon >0$ is arbitrary.
Hence from Proposition \ref{prop:a_mird} we have the desired result.
\end{proof}
We conclude this paper
with a proof of Theorem \ref{thm:mainthm_stnz} (2),
which is much simpler than that of Theorem \ref{thm:h_L_T}.
\begin{thm}
Assume $L_T$ is chosen so that at least one of $L_v$ is a field.
Then $\vartheta_{L_T}^\ird(\mathfrak a,s)$ is holomorphic in the whole complex
plane except for $s=1$ and $5/6$. Also
for any $\varepsilon>0$,
\[
h_{L_T}(\mathfrak a,X)
	=\frac{\gA_{L_T}}{h_k}X
	+\tau(\mathfrak a)\frac{h_k^{(3)}}{h_k}\cdot
		\frac{\gB_{L_T}}{5/6}X^{5/6}
		+O(X^{\frac{4n-1}{4n+1}+\varepsilon})
\qquad (X\to\infty).
\]
\end{thm}
\begin{proof}
The essential fact that we can prove this formula is that
under the condition on $L_T$ we have $\sC(\co,L_T)^\ird=\sC(\co,L_T)$.
By the global and local theory for the space of binary cubic forms
and the orthogonality of characters,
$\vartheta_{L_T}^\ird(\aaa,s)=\vartheta_{L_T}(\aaa,s)$
is holomorphic except for $s=1,5/6$ and
satisfies a functional equation of the form
\[
\vartheta_{L_T}(\aaa,1-s)
=\Gamma(s)^{2n}\Gamma(s-\frac16)^n\Gamma(s+\frac16)^n
\sum_{{\lambda'}\in{\Lambda'}}p_{\lambda'}
(e^{\pi\sqrt{-1}s/2},e^{-\pi\sqrt{-1}s/2})
\xi_{\lambda'}(s)
\]
where $\Lambda'$ is a finite index set and for each $\lambda'\in\Lambda'$,
$p_{\lambda'}(x,y)$ is a polynomial in $x,y$ of degree less than $4n$
and $\xi_{\lambda'}(s)$ is a Dirichlet series with absolute convergence domain
$\re(s)>1$.
Hence the theorem follows from
the modified Landau theorem \cite[Theorem 3]{sash},
the residue formulae of $\vartheta_{L_T}(\aaa,s)$
in Theorem \ref{thm:mainthm_stnz} (1)
and the same treatment for orders of cyclic Galois extensions
as in the proof of Theorem \ref{thm:h_L_T}.
\end{proof}

\end{document}